\documentclass{amsproc}
\usepackage{amssymb}
\usepackage{amsrefs}
\usepackage{bbm}
\usepackage{hyperref}
\newtheorem{theorem}{Theorem}[section]

\newtheorem{proposition}[theorem]{Proposition}

\theoremstyle{definition}
\newtheorem{definition}[theorem]{Definition}
\newtheorem{example}[theorem]{Example}

\newtheorem{corollary}[theorem]{Corollary}
\theoremstyle{remark}
\newtheorem{remark}[theorem]{Remark}
\newtheorem*{sketchproof}{\textit{Sketch of proof}}

\begin{document}

\title[On the LG/CFT correspondence]{On the Landau-Ginzburg/conformal field theory correspondence}

\author[A. Ros Camacho]{Ana Ros Camacho}
\address{Mathematisch Instituut, Universiteit Utrecht\\Postbus 80.010, 3508 TA Utrecht (The Netherlands)}

\email{a.roscamacho@uu.nl}


\date{July 14th, 2018}

\begin{abstract}
In this brief review we introduce the Landau-Ginzburg/conformal field theory correspondence, a result from the physics literature of the late 80s and early 90s which predicts a relation between categories of matrix factorizations and categories of representations of vertex operator algebras. To date there is no precise mathematical statement for it. Here, we explain some recent examples of this correspondence, some standing conjectures and future directions of research.
\end{abstract}

\maketitle

\bibliographystyle{amsplain}

\section{A prediction from physics}

Initially a model to describe superconductivity, Landau--Ginzburg models were promoted in the late 80s to $\left( 2,2 \right)$-supersymmetric quantum field theories completely characterized by a polynomial $W$ called a potential \cite{VW}. These models gained importance in string theory and algebraic geometry as they are a family of quantum field theories which are related under homological mirror symmetry \cites{FJR,Wi}.

Conformal field theories are different kinds of quantum field theories which display conformal symmetry. There has been a substantial effort during the last decades to study this symmetry - rare in nature - and the mathematical structures which encode them. As is well known, this pushed forward our knowledge on (infinite--dimensional) Lie algebras and inspired the definition of vertex operator algebras \cite{Bo} and also that of modular tensor categories.

Despite seeming different, Landau--Ginzburg models and conformal field theories are intimately related. The \textit{Landau-Ginzburg/conformal field theory} (LG/CFT) correspondence is a result from theoretical physics dating from the late 80s and early 90s backed by a generous amount of literature. It states that the infrared fixed point of a Landau--Ginzburg model with potential $W$ is a two--dimensional rational conformal field theory of central charge $c_W$. This implies a direct relation between defects of Landau--Ginzburg models (mathematically described by matrix factorizations) and defects of two--dimensional rational conformal field theories (mathematically described by representations of vertex operator algebras). Hence this correspondence translates to equivalences of categories of matrix factorizations and categories of representations of vertex operator algebras. For suitable subcategories, these equivalences hold not only at the level of $\mathbb{C}$--linear categories but as tensor categories. This is surprising as it connects different areas of mathematics, completely inspired by physics.

There is abundant evidence for these equivalences in several models in conformal field theory \cites{BHLS,KMS,HW1,HW2,HW3,LVW,Mar1,Mar2,VW,Wi}, such as the $N=2$ minimal models or the $N=2$ Kazama--Suzuki model. Mathematically, there exist a promising number of examples \cites{BR,DRCR,CRCR}, with a number of ongoing projects (e.g. \cites{CKRC,RRC,RCW}) providing more in the near future, but we still lack a precise mathematical conjecture for this correspondence.

LG/CFT raises many interesting questions concerning the relation between matrix factorizations and representations of vertex operator algebras. What kind of categorical properties they can share and under which conditions? Which unknown properties can be described on one side suggested by the other? Any advances in this field of research will unveal unknown connections between these two mathematical entities, allowing us to discover common properties not accessible directly.

This brief note is a kind invitation to this fascinating topic for researchers working in the area of vertex operator algebras and their representations -- and so we will assume familiarity of the reader with this topic. Hence, we will first introduce matrix factorizations, and then present some examples of LG/CFT.

\section{Matrix factorizations in a nutshell}

Matrix factorizations were first described by Eisenbud in 1980 in the context of Cohen--Macaulay modules and rings \cite{Ei}. In addition to representations of vertex operator algebras they also have been shown to be related to coherent sheaves \cites{Or1,Or2}, path integrals and quivers \cite{KST}, etc.  In this Section, we give a short introduction to them and the categories one can form with them and their properties.

Let us fix $\mathbbm{k}=\mathbb{C}$, and consider $S:=\mathbbm{k} \left[ x_0,\ldots,x_n \right]$ a polynomial ring and $W \in S$.

\begin{definition}
We call $W$ a \textit{potential} if it satisfies that
$$\mathrm{dim}_{\mathbbm{k}} \left( \frac{S}{\langle \partial_0 W,\ldots,\partial_n W \rangle} \right) < \infty.$$
An equivalent condition is to require that $W$ has an isolated singularity at the origin. 
\end{definition}

\begin{remark}
We constrain ourselves to the particular case where $\mathbbm{k}=\mathbb{C}$ for future convenience when linking to the vertex operator algebras side. It is possible to state the results of this section for more general instances of fields -- at the cost of introducing extra conditions e.g. on $W$ to be a potential. These conditions are, namely, 
\begin{itemize}
\item The quotient ring $\frac{S}{\langle \partial_0 W,\ldots,\partial_n W \rangle}$ is a finitely generated, free $\mathbbm{k}$--module, (which reduces to the above stated one in the case $\mathbbm{k}=\mathbb{C}$),
\item $\partial_{x_i} W$, $i \in \lbrace 0,\ldots,n \rbrace$, form a regular sequence in $S$, and
\item The Koszul complex of $\partial_{x_1} W,\ldots,\partial_{x_n} W$ is exact except in degree zero.
\end{itemize} For more details on the general setting we refer to \cite{CM}.
\end{remark}

\begin{example}
Take $S=\mathbb{C} \left[ x \right]$, and $W=x^d$, for some $d \in \mathbb{Z}_{>0}$. The quotient ring $\frac{\mathbb{C} \left[ x \right]}{\langle x^{d-1} \rangle}$ is clearly finite dimensional as a $\mathbb{C}$--vector space.
\end{example}
\begin{definition}
A \textit{matrix factorization} of a potential $W$ consists of a pair $\left( M,d^M \right)$ where:
\begin{itemize}
\item $M$ is a free $\mathbb{Z}_2$-graded left $S$--module\footnote{$S$ is commutative and hence $M$ can be understood as a left or a right module, but we keep the description of $M$ this way for future convenience when describing a bimodule structure in $M$.}, and
\item $d^M \colon M \to M$ is a degree 1, $S$--linear endomorphism (usually called \textit{twisted differential}) satisfying $d^M \circ d^M=W.\mathrm{id}_M$, where the right hand term stands for the endomorphism $m \mapsto W.m$, for any $m \in M$.
\end{itemize}
We will display the $\mathbb{Z}_2$--grading explicitly as $M=M_0 \oplus M_1$ and $d^M=d_0^M \oplus d_1^M$.
\end{definition}

If $M$ is of finite--rank, then we say that this is a \textit{finite--rank matrix factorization}. Similarly, if $M$ is of infinite--rank it will be called an \textit{infinite--rank matrix factorization}.

\begin{example}
Given a potential $W=x^d$, a non-trivial matrix factorization is the pair $\left( \mathbb{C} \left[ x \right]^{\oplus 2},\begin{pmatrix} 0 & x^m \\ x^{d-m} & 0 \end{pmatrix} \right)$, with $0 < m <d$.
\end{example}

Instead of left modules over $S$, we can also consider bimodules over two polynomial rings $S_1$ and $S_2$ with two potentials $W_1 \in S_1$, $W_2 \in S_2$. In this case
\begin{itemize}
\item the base module is an $S_1$--$S_2$-bimodule instead of a $S$--module\footnote{A $S_1$--$S_2$-bimodule is called free if it is free as a $S_1 \otimes_{\mathbbm{k}} S_2$--module.}, and
\item the twisted differential is $S_1$--$S_2$--bilinear, and it squares to $W_1.\mathrm{id}_M-\mathrm{id}_M.W_2$ (notation: $W_1-W_2$).
\end{itemize}

\begin{example}
Let us take $S_1=\mathbb{C} \left[ x \right]$, $W_1=x^d$, $S_2=\mathbb{C} \left[ y \right]$, $W_2=y^d$. The pair $\displaystyle\left( \mathbb{C} \left[ x \right]^{\oplus 2},\begin{pmatrix} 0 & \prod\limits_{j \in J} \left( x-\eta^{j} y \right) \\ \frac{x^d-y^d}{\prod\limits_{j \in J}  \left( x-\eta^{j} y \right)} & 0 \end{pmatrix} \right)$, where $J \subset \lbrace 0,\ldots,d-1 \rbrace$ and $\eta$ is a primitive $d$-th root of unity, is a matrix factorization of $x^d-y^d$ which we will denote as $P_J$. We call these \textit{permutation-type matrix factorizations}. In the particular case where $J=\lbrace m,\ldots,m+l \rbrace$, denote them as $P_{m:l}$.
\label{permtypeMFs}
\end{example}

It is possible to describe a tensor product of matrix factorizations via this bimodule structure \cite{Yo2}. Let $\left( M,d^M \right)$, $\left( N,d^N \right)$ be two matrix factorizations of $W_1-W_2$ and $W_2-W_3$ respectively. The tensor product matrix factorization $\left( M \otimes_{S_2} N,d^{M \otimes N} \right)$ is a matrix factorization of $W_1-W_3$ where
\begin{itemize}
\item $M \otimes_{S_2} N$ is a $S_1$-$S_3$-bimodule, and
\item $d^{M \otimes N}=d^M \otimes \mathrm{id}_N+\mathrm{id}_M \otimes d^N$.
\end{itemize}

\begin{remark}
Notice that $d^M \otimes \mathrm{id}_N$ and $\mathrm{id}_M \otimes d^N$ are tensor products of graded morphisms. When composing such morphisms, the Koszul sign convention needs to be followed \cite{LV}.
\end{remark}

We distinguish the following categories of matrix factorizations:

\begin{itemize}
\item $\mathrm{MF} \left( W \right)$:
\begin{itemize}
\item[-] Objects are matrix factorizations of the potential $W$, and
\item[-] Given two objects $\left( M,d^M \right), \left( N,d^N \right)$, morphisms of matrix factorizatons consists of all $S$-linear maps $f \colon M \to N$. 
\end{itemize}
This category is clearly differential $\mathbb{Z}_2$-graded. For a morphism of matrix factorizations $f \colon M \to N$ , the differential is given by:
$$\delta f=d^N \circ f- \left( -1 \right)^{|f|} f \circ d^M$$ where $|f|$ is the degree of the morphism $f$.
\item $\mathrm{HMF} \left( W \right)$: 
\begin{itemize}
\item[-] Objects are the same as $\mathrm{MF} \left( W \right)$, and
\item[-] Morphisms are degree zero morphisms of $\mathrm{MF} \left( W \right)$ which lie in the kernel of the differential $\delta$ modulo the image under the differential $\delta$ of the morphisms of degree one.
\end{itemize}
In other words, $\mathrm{HMF}$ is the (zero) homotopy category of $\mathrm{MF}$ \cite{GM}.
\item $\mathrm{hmf} \left( W \right)$: the idempotent complete full subcategory of $\mathrm{HMF} \left( W \right)$ whose objects are (homotopy equivalent to) direct summands of finite--rank matrix factorizations.
\end{itemize}

\begin{remark}
The reason to choose $\mathrm{hmf} \left( W \right)$ to be this particular subcategory is due to the definition of tensor product. Following the notation we used to define the tensor product, when $S_2 \neq \mathbbm{k}$ and $M,N \neq 0$, the tensor product of these two matrix factorizations is infinite--rank. It was proven \cite{DM} that in fact this tensor product is (naturally isomorphic to) a direct summand of some finite--rank matrix factorizations of $W_1-W_3$. Taking the subcategory in this way, we guarantee closure with respect to the tensor product of objects.
\end{remark}

\begin{proposition}[\cites{CR2,CM}]
For the case of matrix factorizations whose base module is a bimodule over $S_1=S_2=S$, $\mathrm{hmf} \left( W \otimes 1-1 \otimes W \right)$ is a tensor category.
\end{proposition}

\begin{example}
In the particular case of $\mathrm{hmf} \left( x^d-y^d \right)$, the unit $I$ of this tensor category is the permutation--type matrix factorization $P_{0:0}$, following the notation of Example \ref{permtypeMFs}. See \cite{CR2} for a detailed description of associator and left and right unitors.
\end{example}

In the case of polynomial rings of one variable \cite{KR,KaR,BR}:
\begin{equation}
\mathrm{Hom}_{\mathrm{hmf} \left( W \otimes 1-1 \otimes W \right)} \left( I,I \right) \cong \frac{S}{\langle \partial_0 W,\ldots,\partial_n W \rangle}.
\nonumber
\end{equation}
A nice discussion on the physical realization of this result, related to the LG/CFT correspondence, can be found at \cite{Va}.

\begin{remark}
\begin{itemize}
\item For categories of matrix factorizations of more general instances of potentials (and defined over more general fields than $\mathbb{C}$) we refer to \cites{CR1,CR2,CM}. In particular, these categories have duals and adjoints whose evaluation and coevaluation maps can be described in a very explicit way, allowing for direct computations. 
\item Categories of matrix factorizations are not in general modular, yet some subcategories (for instance those in the next Section) display interesting structure like semi--simplicity or closedness under tensor product.
\end{itemize}
\end{remark}

Let us conclude this Section by mentioning that a particular composition of the evaluation and coevaluation maps of $\mathrm{hmf}$ is the so--called \textit{quantum dimension}. For a matrix factorization $\left( M,d^M \right)$ of a potential $V \left( x_1,\ldots,x_m \right)-W \left( y_1,\ldots,y_n \right)$, one can compute its left and right quantum dimensions (up to some signs we omit here) via the following formulas:\sloppy
\begin{equation}
\begin{split}
\mathrm{qdim}_l \left( M \right) &= \mathrm{Res} \left( \frac{\mathrm{str} \left( \partial_{x_1} d^M \ldots\partial_{x_m} d^M \partial_{y_1} \ldots \partial_{y_n} \right) dx_1 \ldots dx_m}{\partial_{x_1} V,\ldots,\partial_{x_n} V }\right) \\
\mathrm{qdim}_r \left( M \right) &= \mathrm{Res} \left( \frac{\mathrm{str} \left( \partial_{x_1} d^M \ldots\partial_{x_m} d^M \partial_{y_1} \ldots \partial_{y_n} \right) dy_1 \ldots dy_n}{\partial_{y_1} W,\ldots,\partial_{y_n} W}\right)
\end{split}
\label{qdims}
\end{equation}
where $\mathrm{str}$ is the \textit{supertrace}, defined as the trace of the degree zero part of the resulting matrix minus the trace of the degree one part. Quantum dimensions will be of particular relevance in some of the examples in the next Section.

\section{Examples}

Here we describe in some detail a few available examples of the LG/CFT correspondence.

\subsection{\cite{BR}: the $\mathbf{N=2}$ minimal models}

\quad

For the conformal field theory side, denote by $\nu$ the $N=2$ superconformal vertex operator algebra with central charge $c=3 \left( 1-\frac{2}{d} \right)$, $d \in \mathbb{Z}_{\geq 3}$. Consider the category of representations of the even part $\nu_0$, denoted as $\mathrm{Rep} \left( \nu_0 \right)$. 

\begin{theorem}[\cite{FFRS}]
Consider the Deligne tensor product: $$\mathrm{Rep} \left( \widehat{\mathfrak{su}} \left( 2 \right)_{d-2} \right) \boxtimes \overline{\mathrm{Rep}} \left( \widehat{\mathfrak{u}} \left( 1 \right)_{2d} \right) \boxtimes \mathrm{Rep} \left( \widehat{\mathfrak{u}} \left( 1 \right)_{4} \right)$$ 
where $\mathrm{Rep} \left( \widehat{\mathfrak{su}} \left( 2 \right)_{d-2} \right)$ is the category of integrable highest weight representations of the affine $\mathfrak{su} \left( 2 \right)$ at level $d-2$, $\mathrm{Rep} \left( \widehat{\mathfrak{u}} \left( 1 \right)_{2d} \right)$ is the category of representations for the vertex operator algebra for $\mathfrak{u} \left( 1 \right)$, and the notation $\overline{\phantom{A}}$ means taking the opposite braiding and ribbon twist in this ribbon category. We label simple objects by $\left[ l,m,s \right]$, where $l \in \lbrace 0,\ldots,d-2 \rbrace$, $m \in \mathbb{Z}_{2d}$, $s \in \mathbb{Z}_4$. Then $\mathrm{Rep} \left( \nu_0 \right)$ can be realized as the subcategory with simples such that $l+m+s$ is even\footnote{In other words, we are realizing the category of local modules over a commutative algebra in the triple Deligne tensor product of categories.}. Denote it by $\mathcal{C} \left( d \right)$.
\end{theorem}
$\mathcal{C} \left( d \right)$ has two notable subcategories:
\begin{itemize}
\item $\mathcal{C}_{NS} \left( d \right)$: simples with $l+m$ and $s$ even.
\item $\mathcal{C}_{R} \left( d \right)$: simples with $l+m$ and $s$ odd.
\end{itemize}
The subscripts come from the physical perception of these subcategories, i.e. the Neveu-Schwartz (NS) and the Ramond (R) sectors.

\quad

For the Landau--Ginzburg side, consider the subcategory of $\mathrm{hmf} \left( x^d-y^d \right)$ whose simples are (isomorphic to finite direct sums of) permutation-type matrix factorizations $P_{m:l}$ as defined in Example \ref{permtypeMFs}. Denote it by $\mathcal{PT}_d$.

\begin{theorem}{\cite{BR}}
$\mathcal{C}_{NS} \left( d \right)$ and $\mathcal{PT}_d$ are equivalent as $\mathbb{C}$-linear categories via the assignment: 
\begin{equation}
\left[ l,l+2m,0 \right] \leftrightarrow P_{m:l}
\label{assignment}
\end{equation}
\end{theorem}

\quad

\subsection{\cite{DRCR}: the $\mathbf{N=2}$ minimal models, II}

\quad

The main result of \cite{BR} included a fusion rule, which for matrix factorizations is:
\begin{equation}
P_{m:l} \otimes P_{m':l'} \cong \bigoplus\limits_{k=|l-l'| \quad \mathrm{step} \quad 2}^{\mathrm{min}(l+l',2d-4-l-l')} P_{m+m'-\frac{1}{2} (l+l'-k):k}
\label{fusionrule}
\end{equation}

It is actually possible to make the following improvement:

\begin{theorem}{\cite{DRCR}}
For $d$ odd, $\mathcal{C}_{NS} \left( d \right)$ and $\mathcal{PT}_d$ are equivalent as tensor categories.
\end{theorem}

\begin{sketchproof}
The proof of this result follows several intermediate steps:
\begin{enumerate}
\item For $d$ odd, the category $\mathcal{C}_{NS} \left( d \right)$ is tensor equivalent to the Deligne tensor product of categories $\frac{\mathcal{TL}_d}{\langle p_{d-1} \rangle} \boxtimes \mathrm{Vec}_{\mathbb{Z}_d}$, where $\mathcal{TL}_d$ is the Temperley--Lieb category, $\langle p_{d-1} \rangle$ is the unique ideal tensor generated by the Wenzl--Jones idempotent, and $ \mathrm{Vec}_{\mathbb{Z}_d}$ is the category of $\mathbb{Z}_d$--graded vector spaces.
\item The construction of the tensor functor $F \colon \mathcal{C}_{NS} \left( d \right) \to \mathcal{PT}_d$ amounts to specifying two tensor functors:
\begin{enumerate}
\item a tensor functor $\mathrm{Vec}_{\mathbb{Z}_d} \to \mathcal{PT}_d$,
\item and another $\frac{\mathcal{TL}_d}{\langle p_{d-1} \rangle} \to \left( \mathcal{PT}_d \right)^{\mathbb{Z}_d}$.
\end{enumerate}
Here $\left( \mathcal{PT}_d \right)^{\mathbb{Z}_d}$ is the category of $\mathbb{Z}_d$--equivariant objects in $\mathcal{PT}_d $.
\item In order to construct these functors, it suffices to specify only two pieces of data:
\begin{enumerate}
\item a $\mathbb{Z}_d$--action on $\mathcal{PT}_d$, given by $a \mapsto P_{\lbrace -a \rbrace}$, and
\item  a self--dual object of quantum dimension $2 \cos\left( \frac{\pi}{d} \right)$ in the $\mathbb{Z}_d$--equivariant category of $\mathcal{PT}_d$, which is $P_{\lbrace \frac{d-1}{2},\frac{d+1}{2} \rbrace}$.
\end{enumerate} 
\end{enumerate}
\end{sketchproof}

For $d$ even, this theorem is also expected to hold \cite{RC}. The fusion rule (\ref{fusionrule}) holds for any $d$.

\quad

\subsection{\cite{CRCR}: the $\mathbf{N=2}$ minimal models, III}

\quad

In \cite{CRCR} several equivalences are proven. On the Landau-Ginzburg side, they involve categories of matrix factorizations of potentials describing simple singularities. On the conformal field theory side, categories of modules over algebras representing full conformal field theories in the $N=2$ minimal models. Let us describe these in detail.

\quad

For the conformal field theory side, denote again as $\nu$ the $N=2$ superconformal vertex operator algebra with central charge $c = 3 \left(1 - \frac{2}{d} \right), d \in \mathbb{Z}_{≥3}$, and consider the modular tensor category $\mathrm{Rep} \left( \nu \right)$. It was proven in \cites{FRS,FjFRS} that the full CFTs that can be constructed from a rational vertex operator algebra $\nu$ are parametrized by Morita classes of separable symmetric Frobenius algebras in $\mathrm{Rep} \left( \nu \right)$. For the case of the $N=2$ superconformal vertex operator algebra, we know from \cites{Gan,Gra} that the algebras relevant for LG/CFT are non-trivial only in the $\mathfrak{su} \left( 2 \right)$ factor of the Deligne tensor product category $\mathcal{C} \left( d \right)$. These algebras were identified and classified by Ostrik \cite{Ost}, and are the objects that in the notation of \cite{DRCR} look like:
\begin{equation}
\begin{split}
\left[ 0,0,0 \right] & \oplus \left[ d-2,0,0 \right], \\
\left[ 0,0,0 \right] & \oplus \left[ 6,0,0 \right], \\
\left[ 0,0,0 \right] & \oplus \left[ 8,0,0 \right] \oplus \left[ 16,0,0 \right], \\
\left[ 0,0,0 \right] & \oplus \left[ 10,0,0 \right] \oplus \left[ 18,0,0 \right] \oplus \left[ 28,0,0 \right], \\
\end{split}
\quad \quad
\begin{split}
d &\in 2 \mathbb{Z}_+,\\
d &=12, \\
d &=18, \\
d &=30.
\end{split}
\label{algobj}
\end{equation}

\quad

For the Landau--Ginzburg side, we have the following result.
\begin{theorem}[\cite{CRCR}]
Consider the set of potentials describing simple singularities. These fall into an ADE classification and are as follows:
\begin{equation}
\begin{split}
W_{A_{d-1}} &=x^d+y^2 \quad \quad \quad W_{D_{d+1}}=x^d+x y^2 \\ 
W_{E_6} &=x^3+y^3 \quad \quad W_{E_7} =x^3+x y^3 \quad \quad W_{E_8} =x^3+y^5.
\end{split}
\nonumber
\end{equation}
Then,
\begin{equation}
\begin{split}
\mathrm{hmf} \left( W_{D_{d+1}} \right) & \simeq \mathrm{mod}  \left( P_{\lbrace 0 \rbrace} \oplus P_{\lbrace 0,\ldots,d-1 \rbrace \setminus \lbrace \frac{d}{2}\rbrace}\right) \\
\mathrm{hmf} \left( W_{E_6} \right) & \simeq \mathrm{mod} \left( P_{\lbrace 0 \rbrace} \oplus P_{ \lbrace -3,-2,\ldots,3 \rbrace}\right) \\
\mathrm{hmf} \left( W_{E_7} \right) & \simeq \mathrm{mod} \left( P_{\lbrace 0 \rbrace} \oplus P_{\lbrace -4,-3,\ldots,4 \rbrace} \oplus P_{\lbrace -8,-7,...,8 \rbrace} \right) \\
\mathrm{hmf} \left( W_{E_8} \right) & \simeq \mathrm{mod} \left( P_{\lbrace 0 \rbrace} \oplus P_{\lbrace -5,-4,\ldots,5 \rbrace} \oplus P_{\lbrace -9,-8,\ldots,9 \rbrace} \oplus P_{\lbrace -14,-13,\ldots,14 \rbrace} \right), \\
\end{split}
\nonumber
\end{equation}
where the categories on the right side are categories of modules over the specified algebra objects in the tensor category $\mathrm{hmf}\left( x^d-y^d \right)$.
\label{equivs}
\end{theorem}

Notice that:
\begin{itemize}
\item[-] By construction, the algebras are separable symmetric Frobenius algebra objects in the category $\mathrm{hmf} \left( x^a-y^a \right)$ (where $a=2d,12,18,30$, respectively). 
\item[-] The first equivalence was partly proven in \cite{CR1}.
\item[-] These equivalences are straightforward corollaries of proving an equivalence relation between potentials called \textit{orbifold equivalence}. This equivalence is defined in the context of the bicategory of Landau-Ginzburg models \cites{CR1,CM}. In a nutshell, establishing this equivalence relation between two potentials $V$ and $W$ amounts to finding a finite--rank matrix factorization whose left and right quantum dimensions (as defined in Equation \ref{qdims}) are non--zero. These matrix factorizations have nice categorical properties, yet are computationally challenging to find (as we will see in the next section).
\end{itemize}

\begin{remark}
In the search for these equivalences (and also for most of this note) we work in a setting which includes a $\mathbb{Q}$--grading on the variables of the potential, polynomial ring, and indeed on our matrix factorizations. There are some conditions this grading needs to satisfy, but we have omitted this for the sake of clarity.
\end{remark}

Combining Theorem \ref{equivs} with Assignment (\ref{assignment}) from \cite{DRCR}, we get the following result.

\begin{corollary}
The equivalences of categories of Theorem (\ref{equivs}), under Assignment (\ref{assignment}), read:
\begin{equation}
\begin{split}
\mathrm{hmf} \left( W_{D_{d+1}} \right) & \simeq \mathrm{mod}  \left( \left[ 0,0,0 \right] \oplus \left[ d-2,0,0 \right] \right), \\
\mathrm{hmf} \left( W_{E_6} \right) & \simeq \mathrm{mod} \left( \left[ 0,0,0 \right] \oplus \left[ 6,0,0 \right] \right), \\
\mathrm{hmf} \left( W_{E_7} \right) & \simeq \mathrm{mod} \left( \left[ 0,0,0 \right] \oplus \left[ 8,0,0 \right] \oplus \left[ 16,0,0 \right] \right), \\
\mathrm{hmf} \left( W_{E_8} \right) & \simeq \mathrm{mod} \left( \left[ 0,0,0 \right] \oplus \left[ 10,0,0 \right] \oplus \left[ 18,0,0 \right] \oplus \left[ 28,0,0 \right] \right), \\
\end{split}
\quad \quad
\begin{split}
d &\in 2 \mathbb{Z}_+,\\
d &=12, \\
d &=18, \mathrm{and} \\
d &=30.
\end{split}
\nonumber
\end{equation}
\label{orbcor}
The algebra objects match exactly those from Equation \ref{algobj}.
\end{corollary}

\quad

\subsection{Beyond the $N=2$ minimal models}

\quad

Proceeding analogously to the previous Section, one may wonder if this equivalence relation between different potentials could provide further examples of LG/CFT. Thanks to Arnold, we have a complete classification of simple, unimodal and bimodal singularities \cite{Ar,AGV}, some of which may serve as a list of potential candidates for orbifold equivalence (provided some conditions are satisfied). Physics literature \cite{Mar1,Mar2} suggests that examples involving potentials describing singularities beyond the simple case could be related to more exotic models than the minimal ones previously discussed. Following this lead,  we list a collection of results recently achieved in this direction.

\begin{itemize}

\item[\cite{RCN1,RCN2}:] In the first paper, two equivalences between potentials describing exceptional unimodal singularities were proven. These two singularities are related by the so--called Arnold strange duality \cite{Ar,Eb}. Similarly, in \cite{RCN2} equivalence between different potentials describing the same exceptional unimodal singularities were proven. In each of these cases, we find equivalences depending on parameters satisfying a list of equations whose solutions are related by Galois groups (which were also discussed in \cite{DRCR} and \cite{CRCR}). Galois groups were described in rational CFTs in \cite{Gep1,Gep2}. 
\item[\cite{KCL}:] The remaining instances of equivalences involving potentials describing exceptional unimodal singularities related by strange duality were proven, as well as several other instances of potentials including some of Fermat, chain and loop type (for more information on these types, see \cite{KS,Her}). A discussion on the computational challenges posed by this search of equivalences is included.
\item[\cite{CKRC}:] This work proves that the computational difficulties found in \cite{CRCR,RCN1,RCN2} and discussed in \cite{KCL} are not by chance. A systematic algorithm analysis, used to find the equivalences, shows one needs to solve an amount of variables and equations way beyond our current computational powers. Nevertheless, a ``lucky guessing'' strategy -- simplifying these amounts -- can be performed, giving satisfactory results. This is shown for two extra equivalences involving potentials describing exceptional unimodal and bimodal singularities.
\end{itemize}

There is current work in progress \cite{RC} which aims to provide a CFT interpretation of these results in the same way as Corollary \ref{orbcor}.

\quad

\subsection{Further examples}

\quad

We would like to conclude our exposition mentioning some works in progress which in the very near future will potentially shed some light on a precise mathematical statement for LG/CFT.
\begin{itemize}
\item[\cite{RRC}:] This project attacks the problem of extending the main theorem of \cite{DRCR} to the subcategory $\mathcal{C}_R \left( d \right)$. The way to do this is via giving a precise description of spectral flow, which one can find comparing orbits of objects on each category. A precise description of the conjugation automorphisms present in the category of matrix factorizations is also included.
\item[\cite{KRC}:] It is a well--known fact that characters of representations of vertex operator algebras are modular forms. In the setting of the $N=2$ minimal models, they have been known for a long time but have been recently corrected \cite{CLRW}: they are described in terms of Jacobi theta and Dedekind eta functions. The aim of this project is to describe the analogous concept to these characters for matrix factorizations.
\item[\cite{RCS}:] The previous works are mostly focused on the $N=2$ minimal models. There is physical evidence of similar results in $N=2$ Kazama--Suzuki models \cites{BF1,BF2}, which should provide insights to the dependence of LG/CFT on the model we are considering. This project gives a precise mathematical description of these physical results.
\item[\cite{RCW}:] The target of this project is to provide a functorial description of LG/CFT between suitable sub(bi)categories which should include all of the previous results. This would be the summit of all of the aforementioned projects. The starting point is the construction of a higher category whose objects are pairs of a (rational, $c_2$--cofinite) vertex operator algebra and a Frobenius algebra object in its category of representations. This construction is similar to the 2--category of Frobenius algebras described in \cite{CSP}, and is inspired by \cite{FFRS,FRS}. Although we have omitted it in this note, it is possible to formulate the Landau--Ginzburg side in the shape of another 2--category whose objects are pairs of polynomial rings and potentials in them (and actually, reformulate all the results mentioned on matrix factorizations in an elegant way). With these two, the aim is to describe a higher functor which should be the precise mathematical formulation of the LG/CFT correspondence.
\end{itemize}

\end{document}